\documentclass[12pt,twoside,reqno,psamsfonts]{amsart}
\usepackage{amssymb}
\makeindex
\def\R{{\mathbb R}}

\def\N{{\mathbb N}}

\def\squareforqed{\hbox{\rlap{$\sqcap$}$\sqcup$}}
\def\qed{\ifmmode\squareforqed\else{\unskip\nobreak\hfil
\penalty50\hskip1em\null\nobreak\hfil\squareforqed
\parfillskip=0pt\finalhyphendemerits=0\endgraf}\fi}

\begin{document}

\title{Slice regular Malmquist-Takenaka system in the quaternionic Hardy spaces}
\author{\sc Margit Pap}
\address{\textnormal{Faculty of Sciences, University of P\'ecs, 7634 P\'ecs, Ifj\'us\'ag \'ut 6, HUNGARY} \newline
\textnormal{e-mail: papm@gamma.ttk.pte.hu}}
\begin{abstract} A slice regular  analogue of the Malmquist-Takenaka system  is investigated.  
It is proved that they form a complete orthonormal  system in the quaternionic Hardy spaces of the unit ball.  The properties of associated projection operator are studied.

{\it Key Words and Phrases}: Functions of hypercomplex variables, quaternionic Hardy spaces, series expansions, interpolation, approximation by rational functions, quaternionic Malmquist-Takenaka system

{\it Mathematics Subject Classification (2010)}: 30G35, 30H10, 30D55,41A05,   41A20, 41A58
\end{abstract}
\maketitle


\section{ Introduction}
\vskip5mm

The first mention of rational orthonormal systems  in the Hardy space of complex variable functions seems to have occurred  in the work of Takenaka and Malmquist \cite{M, T}. These systems can be viewed as  extensions of the trigonometric system on the unit circle, that corresponds to the special choice when  all of the poles are located at the origin. In the system theory  they are often used to identify the transfer function of the system.

This orthonormal system is generated by a sequence $a=(a_1,a_2,...)$ of complex numbers, $a_n \in \mathbb{D}$ of the unit disc $\mathbb{D} :=\{ z\in \mathbb{C} : |z|<1 \}$ and can be expressed by the Blaschke-functions
$$
B_b(z):=\frac{z-b}{1-\overline{b}z} \quad \quad (b\in \mathbb{D} ,\,  z\in \mathbb{C}).
$$
If $b$ belongs to $\mathbb{D}$, then $B_b$ is a one-to-one map on $\mathbb{D}$ and $\mathbb{T}$, respectively.

The Malmquist-Takenaka system (M-T)  $\Phi_n=\Phi_n^a \,(n\in \N^*) $   is defined by
$$
\Phi_1(z)=\frac{\sqrt{1-|a_1|^2}}{1-\overline{a_1}z},
$$
$$
\Phi_n(z)=\frac{\sqrt{1-|a_n|^2}}{1-\overline{a_n}z}\prod_{k=1}^{n-1}B_{a_k}(z),\,\, n\ge 2.
$$
These  functions form an orthonormal system on the unit circle $\mathbb{T}:=\{ z\in \mathbb{C} : |z|=1\}$, i.e., 
$$
\langle \Phi_n, \Phi_m \rangle=\frac1{2\pi}\int_0^{2\pi}\Phi_n(e^{it}) \overline{\Phi_m(e^{it})}dt=\delta_{mn} \quad (m,n\in \mathbb{N}^*),
$$
where $\delta_{nm}$ is the Kronecker symbol.
If the sequence  $a=(a_1,a_2,...)$ satisfies the non-Blaschke condition
$$
\sum_{n\ge 1}(1 - |a_n|) = +\infty,$$
then the corresponding M-T system is complete in the Hardy space of the unit disc.

In the papers \cite{pasc01, pa03} there was proved a discrete orthogonality property of these functions, and a quadrature method based on interpolation formula was given. 
When all the parameters are equal, i.e., $a_n=a, \, n\in \mathbb{N}^*$, we obtain the so called discrete Laguerre system and particularly when $a_n=0, \, n\in \mathbb{N}^*$ we obtain the trigonometric system.

In the  papers \cite{pa04},  \cite{qian12} M. Pap, F.Schipp and Qian T, Sprossig W, Wang J. respectively, following two different ways,   introduced two analogues of the M-T systems in the set of quaternions. The drawback of both constructions is that these extensions will not  inherit all the nice properties of the before mentioned system, e.g., the system introduced by Pap and Schipp is not analytic in the quaternionic setting. The system introduced by Qian T, Sprossig W, Wang J., is monogenic   but  can not  be  written in closed form.

In this paper we introduce the slice regular analogue of the M-T system and we will study its properties.
 
\vskip5mm
\section{ Quaternions}
\vskip5mm

 Quaternions are  extensions of complex numbers. Instead of just $i$ we have three numbers that all are square roots of $-1$, denoted by $i,j$ and $k$ satisfying the following identities:
$$
i^2=j^2=k^2=-1,\quad ij=-ji=k, \quad ki=-ik=j,
$$
named as Hamilton's rules.
A quaternion $q$ can be represented as 
$$
q=x_1+x_2i+x_3j+x_4k, \quad (x_n \in \R,\, n=1,2,3,4).
$$
Let the set of quaternions be denoted by $${\mathbb{H}}:=\{q=x_1+x_2i+x_3j+x_4k: \quad x_n \in \R,\, n=1,2,3,4\}.$$

The sum of two quaternions is:
$$
a+b=(a_1+b_1)+(a_2+b_2)i+(a_3+b_3)j+(a_4+b_4)k.
$$
From  Hamilton's rules it follows that the product of two quaternions is 
$$
a.b=(a_1b_1-a_2b_2-a_3b_3-a_4b_4)+(a_1b_2+a_2b_1+a_3b_4-a_4b_3)i+ $$ $$ 
(a_1b_3-a_2b_4+a_3b_1+a_4b_2)j+(a_1b_4+a_2b_3-a_3b_2+a_4b_1)k.
$$

The conjugate of a quaternion $q$  is given by 
$$
\overline{q}=x_1-x_2i-x_3j-x_4k,
$$
and the quaternion norm is
$$
\|q\|=\sqrt{q.\overline{q}}=\sqrt{\overline{q}.q}=\sqrt{x_1^2+x_2^2+x_3^2+x_4^2}.
$$
We mention that the product of two quaternions is not commutative in general and $\overline{a.b}=\overline{b}.\overline{a}$.
The multiplicative inverse of $q$ is $$ q^{-1}=\frac{\overline{q}}{q.\overline{q}}.$$
$(\mathbb{H}, +,.)$ being a noncommutative field (skew field). 

In these notations it can be seen that  quaternions are closely related to four dimensional
vector spaces, moreover as $\mathbb{C}$ can be identified with $\mathbb{R}^2$, in the same manner $\mathbb{H}$ can be identified with $\mathbb{R}^4$. 

While in $\mathbb{C}$ the equation $z^2+1=0$ has two solutions $i$ and $-i$, in $\mathbb{H}$ this equation has infinitely many quaternion solutions. Namely all ${ I}=u_2i+u_3j+u_4k$, with $u_2^2+u_3^2+u_4^2=1$,
satisfy this equation. Geometrically we can interpret this in the following way: every three dimensional unit vector in the system $(i,j,k)$ is a solution of $z^2+1=0$.

It is well known that  every unit complex  number $z$ can be represented in the trigonometric form as :
$$
z=\cos \theta +i \sin \theta=e^{i\theta}\quad (\theta \in [0,2\pi)).
$$
The analogue of this formula holds for quaternions.
Indeed set 
$$
e^{q}:=\sum_{n=0}^{\infty}\frac{q^n}{n!}\qquad (q\in \mathbb{H}).
$$
If $I$ is a unit quaternion of the form ${ I}=u_2i+u_3j+u_4k$ with $u_2^2+u_3^2+u_4^2=1$, then ${ I}^2=-1$. It follows that $${I}^{4l}=1, \quad {I}^{4l+1}={I}, \quad {I}^{4l+2}=-1, \quad {I}^{4l+3}=-{I} \quad (l\in \mathbb{N}).$$

Hence we obtain that for every $\theta_1 \in [0, 2\pi)$
$$
e^{\theta_1 {I}}=\sum_{n=0}^{\infty}\frac{(\theta_1 I)^n}{n!}=\sum_{n=0}^{\infty}\frac{(-1)^n \theta_1^{2n}}{(2n)!}+\sum_{n=0}^{\infty}\frac{(-1)^n\theta_1^{2n+1}}{(2n+1)!}{I}
$$
$$
=\cos \theta_1 +I .\sin \theta_1 .
$$
Using  spherical coordinates ${I}$ can be expressed as 
$$
{I}=\cos \theta_2i+\sin \theta_2 \cos \theta_3 j +\sin \theta_2 \sin \theta_3k,\quad \theta_3 \in [0, 2\pi),  \theta_2 \in [0, \pi]. 
$$

Let  us denote by $\mathbb{B}$ the open unit ball in the set of quaternions,  its boundary by $\partial \mathbb{B}$. Let also denote by $$\partial{\mathbb{B}_I}:=$$ $$\{ q\in \mathbb{H}: q=\cos \theta_1+{I}.\sin \theta_1, \, {I}=\cos \theta_2i+\sin \theta_2 \cos \theta_3 j +\sin \theta_2 \sin \theta_3k,$$
$$
 \theta_3 \in [0, 2\pi), \theta_1, \theta_2 \in [0, \pi]\}$$  the set of unit quaternions defined by the direction $I$. $(\partial \mathbb{B}_I, .)$ is a compact commutative  subgroup of $(\mathbb{H}, .)$.

The Haar-measure on this subgroup is defined by
$$
\int_{\partial \mathbb{B}_I}f(q)dq:=\frac1{\omega_3}\int_0^{2\pi}\int_0^{\pi}\int_0^{\pi}f(q)\sin^2 \theta_1\sin \theta_2 d\theta_1 d\theta_2 d\theta_3,
$$
where 
$$
\omega_3=\int_0^{2\pi}\int_0^{\pi}\int_0^{\pi}\sin^2 \theta_1\sin \theta_2 d\theta_1 d\theta_2 d\theta_3=4\pi^2.
$$

In engineering applications  unit quaternions are used to describe
 the rotations. Namely the unit quaternion $q$ corresponding to a rotation around vector $$I=(\cos \theta_2,\sin \theta_2 \cos \theta_3 ,\sin \theta_2 \sin \theta_3)$$ by angle $2\theta_1$ becomes  $q=\cos \theta_1+{I}.\sin \theta_1$. If $X=x_1i+x_2j+x_3k$ is the corresponding quaternion of the vector $x=(x_1,x_2,x_3)$ then the transformation $$ Y=qX\overline{q}$$ is a rotation around vector $I$ by angle $2\theta_1$.

Another important field, where  quaternions play important role is the quantum theory.

\vskip5mm

\section{ Slice regular functions}
\vskip5mm

One way to generalize  Blaschke functions over the set of  quaternions would be the direct extension of the before mentioned formula over the quaternions, i.e., 
\begin{equation}
\label{Bclas}
B_a(q):=(1-q\overline{a})^{-1}{(q-a)}=\frac1{\lambda_{a,q}}\overline{(1-q\overline{a})}(q-a)=\frac1{\lambda_{a,q}}(1-a\overline{q})(q-a)
\end{equation}
where 
$$
\lambda_{a,q}=\overline{(1-q\overline{a})}(1-q\overline{a})={(1-a\overline{q})}(1-q\overline{a}), a\in \mathbb{B}, q\in \overline{\mathbb{B}}.
$$

Recently in \cite{alco12, bist13, st10} it was introduced and studied the so called slice regular Blaschke functions.

The theory of slice regular functions of a quaternionic variable (often simply called regular functions) was
introduced in  \cite{gest06}, \cite{gest07}, and represents a natural quaternionic counterpart of the theory of complex holomorphic
functions. This recent theory has been growing very fast and was developed in a series of papers,
including in particular \cite{coge10}, \cite{coge09}, \cite{cosa11}, \cite{cu65}, \cite{gest12}, \cite{st10} where most of the recent
advances are discussed. The detailed up-to-date theory appears in the monograph \cite{gest13}.
The theory of regular functions is presently expanding in many directions.

 Set $\mathbb{S} = \{q \in \mathbb{H} : q^2 = -1\}$ to be the $2$-sphere of purely
imaginary units in $\mathbb{H}$, and for $I \in \mathbb{S}$ let $L_I$ be the complex plane $\mathbb{R} + \mathbb{R}I$, then we have 
$$\mathbb{H} =\cup_{I\in \mathbb{S}}L_I .$$

To recall the definition of slice regular function we will first describe the natural domains of definition for
such functions (for the definitions and main results see the monograph  \cite{gest13} and the reference list therein).

{\bf Definition 3.1.} Let $\Omega$
 be a domain in $\mathbb{H}$ that intersects the real axis. Then:
 
 1. $\Omega$ is called a slice domain if, for all $I \in  \mathbb{S}$, the intersection 
$\Omega_I$ with the complex plane $L_I$ is a domain
of $L_I$ ;
 
2. $\Omega$ is called a symmetric domain if for all $x, y \in  \mathbb{R}, x + yI \in \Omega$
 implies $x + y\mathbb{S} \subset \Omega$.

{\bf Definition 3.2.} Let 
$\Omega \subset \mathbb{H}$ be a slice domain. A function $f : \Omega \to \mathbb{H}$ is said to be (slice) regular if, for all
$I \in \mathbb{S}$, its restriction $f_I$ to 
$\Omega_I$ is holomorphic, i.e., it has continuous partial derivatives and satisfies
\begin{equation}\overline{\partial}_If(x + yI) :=\frac{1}{2} \left(\frac{\partial}{\partial x}+ I \frac{\partial}{\partial y}\right)f_I (x + yI) = 0.
\end{equation}

{\bf Lemma 3.1} (Splitting Lemma). If f is a regular function on a slice domain 
, then for every $I \in \mathbb{S}$ and for
every $J \in \mathbb{S}$, $J$ orthogonal to $I$, there exist two holomorphic functions $F,G : \Omega_I \to L_I$ , such that for every
$z = x + yI \in I$, we have
\begin{equation}
f_I (z) = F(z) + G(z)J.
\end{equation}

As shown in \cite{gest07}, if 
 we consider the open unit ball $\mathbb{B}$ of $\mathbb{H}$, the class of regular functions coincides with the class of
convergent power series of type
$\sum_{n\ge 0} q^na_n$, with all $a_n \in\mathbb{H}$.

The direct extension  of the Blaschke function, presented before, is not slice regular. In general the product of two slice regular functions is not slice regular.

{\bf Definition 3.3.} Let $f, g :\mathbb{B}\to \mathbb{H}$ be regular functions and let $f(q)=\sum_{n\in \mathbb{N}} q^na_n$; $g(q)=\sum_{n\in \mathbb{N}} q^nb_n$ be their power series expansions. The regular
product of $f$ and $g$ (sometimes referred to as their $\ast$-product) is the regular function
defined by
\begin{equation}
f \ast g(q)=\sum_{n\in \mathbb{N}}q^n \sum_{k=0}^na_kb_{n-k} 
\end{equation}
on the same ball $\mathbb{B}$.

We can define two additional operations on regular functions. 

{\bf Definition 3.4.} Let $f:\mathbb{B}\to \mathbb{H}$ be a regular function and let $f(q)=\sum_{n\in \mathbb{N}} q^n a_n$ be its power series expansion. The regular conjugate of $f$ is the regular
function defined by
$f^c(q)=\sum_{n\in \mathbb{N}}q^n \overline{a_n}$ 
on the same ball $\mathbb{B}$. The symmetrization of $f$ is the function
$f^s = f\ast f^c= f^c\ast f$.

{\bf Definition 3.5.} Let $f$ be a regular function on a symmetric slice domain $\Omega$.
If $f \neq 0$ on $\Omega$, the regular reciprocal of $f$ is the function  
$$f^{-\ast}=(f^s)^{-1}f^c. $$

{\bf Definition 3.6.} The  regular Blaschke function by definition is:

\begin{equation}
\mathcal{B}_a(q) = (1 - q \overline{a})^{-\ast}\ast (q - a).
\end{equation}

This function inherits all the nice properties of the complex Blaschke functions, i.e., is a regular fractional transformations that maps the open quaternionic unit ball $\mathbb{B}$
onto itself and the boundary of unit ball $\partial \mathbb{B}$
onto itself bijectively (see \cite{alco12, bist13, st10}).

The  classical and regular Blaschke functions are related in the following way:
$$\mathcal{B}_a(q) = (1 - q \overline{a})^{-\ast}\ast (q - a)=B_a(T_a(q)),$$
where $T_a (q) = (1 - qa)^{-1}q(1 - qa)$ is a diffeomorphism of $\mathbb{B}$.

It can also be proved  that the factors in the definition of the regular Blaschke product commute
$$
\mathcal{B}_a(q) =(1 - q \overline{a})^{-\ast}\ast (q - a)= (q - a)\ast(1 - q \overline{a})^{-\ast}.$$

One of the most fertile chapters of the theory of complex holomorphic functions consists of the theory of
Hardy spaces. This theory contains results of great significance that led to important general achievements
and subtle applications, i.e., the transfer function of a linear time invariant system belongs to the Hardy space. In the  paper \cite{fage14}   the quaternionic counterpart of complex Hardy spaces was introduced, and  their basic and fundamental properties was 
investigated. 

{\bf Definition 3.7.} Let $f : \mathbb{B} \to \mathbb{H}$ be a regular function and let $0 < p < +\infty$. Set
\begin{equation}||f||_p = \sup_{I\in \mathbb{S}}\lim_{r\to 1_-}\frac{1}{2\pi}\left(\int^{2\pi}_0|f(re^{I\theta})|^pd\theta\right)^{1/p},
\end{equation}
and set
\begin{equation}
||f||_{\infty} = \sup_{q\in\mathbb{B}}|f(q)|.
\end{equation}
Then, for any $0 < p \geq +\infty$, we define the quaternionic Hardy space $H^p(\mathbb{B})$ as
\begin{equation}
H^p(\mathbb{B}) = \{f : \mathbb{B}\to  \mathbb{H} | f\,  \text{is\, \, regular\,\, and} \,\,||f||p < +\infty\}.
\end{equation}
In \cite{fage14}  the main properties and features of the quaternionic $H^p$-norms were  studied,  and  the initial properties of the quaternionic $H^p$ spaces were established. It turned out that many properties of the complex Hardy spaces have their quaternionic analogue. The
boundary behavior of functions $f$ in $H^p(\mathbb{B})$  is very similar to the complex case, i.e.,  for almost every $\theta \in  [0, 2\pi)$, the limit
\begin{equation}
\lim_{r\to 1_-}f(re^{I\theta}) =  \widetilde{f}_I (e^{I\theta})
\end{equation}
exists for all $I \in \mathbb{S}$ and in this case it belongs to $L^p(\partial B_I )$. The properties of the boundary
values of the $\ast$-product of two functions, each belonging to some $\mathbb{H}^p(\mathbb{B})$ space was  investigated. The
classical $H^p$ kernels, and in particular the Poisson kernel, to the quaternionic setting was  extended  and  the Poisson-type and Cauchy-type 
representation formulas for all $f \in  H^p(\mathbb{B})$ was deduced. Analogues of outer and inner functions
and singular factors on $\mathbb{B}$ were given, whose definitions (when compared with those used in the complex case)
clearly resent of the peculiarities of the non commutative quaternionic setting. 
Factorization properties of $H^p$ functions were established. The Blaschke factor of a function $f$ in
$H^p(\mathbb{B})$ are built from the zero set of $f$ using the regular Blaschke functions; it can be obtained also a complete factorization result, in terms of an outer,
a singular and a Blaschke factor, for a subclass of regular functions, namely for the one-slice-preserving
functions.

It is valid the following Splitting Formula: if $f \in  H^p(\mathbb{B})$ for some $p \in  (0,+\infty]$, then for any $I \in \mathbb{S}$,  the splitting of $f$ on $L_I$
with respect to $J \in \mathbb{S}$, $J$ orthogonal to $I$, is $f_I (z) = F(z) + G(z)J$, then the holomorphic functions $F$ and $G$ are both
in $H^p(\mathbb{B}_I )$.

In analogy with the complex case, the space $H^2(\mathbb{B})$ is special. Indeed the $2$-norm turns out to be induced
by an inner product.
 Let $f \in H^2(\mathbb{B})$ and let $f(q) =\sum_{n\ge 0}q^na_n$ be its power series expansion. Then the
square of the $2$-norm of $f$,
coincides with 
\begin{equation}
\|f\|^2=\sum_{n\ge 0}|a_n|^2.
\end{equation}

This result  permits a  way to define an inner product on the space $H^2(\mathbb{B})$. In fact, if $f, g \in H^2(\mathbb{B})$,
let $f(q) =\sum_{n\ge 0} q^na_n,\, g(q) =\sum_{n\ge 0} q^nb_n$ be their power series expansions, based on  previous result, then their inner product is defined by
\begin{equation}
\langle f, g\rangle=\lim_{r\to 1}\frac{1}{2\pi}\int_0^{2\pi}\overline{g(re^{I\theta})}f(re^{I\theta})d\theta =\sum_{n \ge 0}\overline{b_n}a_n,
\end{equation}
 for any $I \in \mathbb{S}$.

Thanks to the existence of the radial limit it is possible to obtain integral representations for functions in
$H^p(\mathbb{B})$ for $p \in [1,+\infty]$.

{\bf Theorem 3.1.} If $f \in  H^p(\mathbb{B})$ for $p \in [1,+\infty]$, then, for any $I \in \mathbb{S}$, $f_I$ is the Poisson integral and the Cauchy
integral of its radial limit $\tilde {f}_I$ , i.e.,
\begin{equation}
f_I (re^{I\theta}) =\frac{1}{2\pi}\int_0^{2\pi}\frac{1 - r^2}{1 - 2r cos(\theta - t) + r^2}\tilde{f}_I(e^{It})dt
\end{equation}
and
\begin{equation}
f_I (z) =\frac{1}{2\pi I}\int_{\partial B_I}\frac{d\xi}{\xi-z}\tilde{f}_I(\xi).
\end{equation}

The next result, on the other hand, is a more powerful Cauchy Formula, which
allows the reconstruction  of $f$ on the entire open set of definition, by
using its values on a  single slice.

{\bf Theorem 3.2} (Cauchy Formula). Let $f$ be a regular function on a symmetric
slice domain $\Omega$. If $U$ is a bounded symmetric open set with $U \subset \Omega$,  $I \in \mathbb{S}$, and if
$\partial U_I$ is a finite union of disjoint rectificable Jordan curves, then, for $q \in  U$,
\begin{equation}
f(q)=\frac{1}{2\pi}\int_{\partial U_I}(s-q)^{-\ast}ds_If_I(s),
\end{equation}
where $ds_I=-I ds$ and $(s-q)^{-\ast}$ denotes the regular reciprocal of $(s-q)$.

{\bf Theorem 3.3.} (Zero set structure). Let $f$ be a regular function on a symmetric slice domain. If $f$ does
not vanish identically, then its zero set consists of the union of isolated points and isolated $2$-spheres of the
form $x + y\mathbb{S}$ with $x, y \in \mathbb{R}, y \ne 0$

Spheres of zeros of real dimension $2$ are a peculiarity of regular functions.

Let $f$ be a regular function on a symmetric slice domain. A $2$-dimensional sphere
$x + y\mathbb{S} \subset {\mathcal{Z}}_f$ 
 of zeros of $f$ is called a spherical zero of $f$ and is represented by an element $x + yI$ of such a sphere, called a
generator of the spherical zero $x + y\mathbb{S}$. Any zero of $f$ that is not a generator of a spherical zero is called an
isolated zero (or a non spherical zero or simply a zero) of $f$.

{\bf Theorem 3.4.} Let $p \in (0,+\infty], f \in H^p(\mathbb{B}), f \ne 0$ and let $\{b_n\}_{n\in\mathbb{N}}$ be its sequence of zeros. Then
$\{b_n\}_{n\in\mathbb{N}}$ satisfies the Blaschke condition
$$
\sum_{n\ge 0}(1 - |b_n|) < +\infty.$$

Holomorphic functions defined on a domain 
$\Omega_I$, symmetric with respect to the real axis in the complex plane
$L_I$, extend uniquely to the smallest symmetric slice domain of $\mathbb{H}$ containing 
$\Omega_I$.

{\bf Theorem 3.5.} (Extension Lemma). Let $\Omega$
 be a symmetric slice domain and choose $I \in \mathbb{S}$. If $f_I : \Omega_I \to \mathbb{H}$ is
holomorphic, then setting
\begin{equation}
f(x + yJ) =\frac{1}{2}[f_I (x + yI) + f_I (x - yI)] + \frac{JI}{2}[f_I(x - yI) - f_I (x + yI)]
\end{equation}
extends $f_I$ to a regular function $f : \Omega \to \mathbb{H}$. Moreover $f$ is the unique extension and it is denoted by $ext(f_I )$.


\section{ The regular quaternionic  Malmquist-Takenaka system  }
\vskip5mm

 The extension of the M-T systems for quaternions in  \cite{pa04} was  described as follows: let us consider  a sequence $a=(a_1,a_2,...)$ of quaternions, $|a_n|<1,\,(n\in \N^*)$  and the classical quaternionic extension of Blaschke-functions $B_{a}$.
The functions  $\Phi_n=\Phi_n^a \,(n\in \N^*) $  are defined very similar to the complex case by the quaternionic product 
 
$$
\Phi_1(z)={\sqrt{1-|a_1|^2}}(1-z\overline{a_1})^{-1}, 
$$
\begin{equation}
\Phi_n(z)={\sqrt{1-|a_n|^2}}\left (\prod_{k=1}^{n-1}B_{a_k}(z)\right )(1-z\overline{a_n})^{-1}\,\,\,(z\in \overline {\mathbb{B} }, \, n=2,3,...).
\end{equation}
Unfortunately, they are not regular functions anymore. But still  for their Dirichlet kernel it was possible to prove the analogue of the Darboux-Christoffel formula.

When all the parameters are equal $a_n=a=re^{\theta {\bf I}}=r(\cos \theta+{\bf I}\sin \theta) \, (n\in \N^*)$, we obtain the quaternionic analogue of the discrete Laguerre system. Even in this special case the orthogonality is not yet proved, but  a discrete orthogonality property of this particular case can be proved(see \cite{pa04}).

In \cite{qian12}  the authors studied the decompositions of functions in the quaternionic monogenic  Hardy spaces into linear combinations of the basic functions in the orthogonal rational systems,  which can be obtained in the respective contexts through Gram-Schmidt orthogonalization process on shifted
Cauchy kernels.  While in the complex case, following these two ways we get the same system, here in the quaternionic case it has not been  proved yet,   that the two methods  give  the same.

In this paper we will consider the slice regular analog of the Malmquist-Takenaka system and we will investigate the properties of this system.

Let consider  a sequence $a=(a_1,a_2,...)$ of quaternions, $|a_n|<1,\,(n\in \mathbb{N}^*)$.  The slice regular analogue of the Malmquist-Takenaka system can be expressed by the slice regular quaternionic Blaschke-functions $\mathcal{B}_{a_n}(q)=(1 - q \overline{a_n})^{-\ast}\ast (q - a_n)= (q - a_n)\ast(1 - q \overline{a_n})^{-\ast}$.
Namely, the functions  $\bf{\Phi}_n=\bf{\Phi}_n^a \,(n\in \N^*) $  are defined very similar to the complex case by, but here we use the slice regular product of the factors:
$$
{\bf{\Phi}}_1(z)={\sqrt{1-|a_1|^2}}(1-z\overline{a_1})^{-\ast},
$$
\begin{equation}
{\bf{\Phi}}_n(z)={\sqrt{1-|a_n|^2}}\left (\ast\prod _{k=1}^{n-1}{\mathcal{B}}_{a_k}(z)\right )\ast (1-z\overline{a_n})^{-\ast}\,\,\,(z\in \overline{\mathbb{B}}, \, n=2,3,...),
\end{equation}
where  $\ast\prod$ means the $\ast$-product of the factors. Because $\mathcal{B}_{a}(q)$ is a slice  regular function and the  $\ast$-product of two slice regular functions is slice regular, in this way we generate a slice regular system.

 When all the parameters are equal $a_n=a=re^{\theta {\bf I}}=r(\cos \theta+{\bf I}\sin \theta) \, (n\in \N^*)$, then we get the slice regular analogue of the discrete Laguerre system, 
$$
{\bf{L}}_1(z)={\sqrt{1-|a|^2}}(1-z\overline{a})^{-\ast},
$$
$$
{\bf{L}}_n(z)={\sqrt{1-|a|^2}}\left (\ast\prod _{k=1}^{n-1}{\mathcal{B}}_{a}(z)\right )\ast (1-z\overline{a})^{-\ast}\,\,\,(z\in {\overline{\mathbb{B}}}, \, n=2,3,...),
$$
When all the parameters are $0$, i.e.  $a_n=0 \, (n\in \N^*)$ we obtain ${\bf{\Phi}}_n(z)=z^n$,   the quaternionic analogue of the trigonometric system.

{\bf Lemma 4.1.}
The slice regular analogue of the discrete Laguerre system can be written in the following form:
$$
{\bf{L}}_n(z)={\sqrt{1-|a|^2}}(q - a)^{\ast n}\ast(1 - q \overline{a})^{-\ast (n+1)} (z\in \overline{\mathbb{B}}, \, n=2,3,...).
$$
 
This can be proved by induction, using the commutativity property of the factors in   $\mathcal{B}_{a}(q)=(1 - q \overline{a})^{-\ast}\ast (q - a)= (q - a)\ast(1 - q \overline{a})^{-\ast}$.

{\bf Theorem 4.1.} 
If all the parameters of the slice regular Malmquist -Takenaka system are on the same slice, i.e., there exists $I\in \mathbb{S}$ such that $a_n=r_ne^{\theta_n { I}}=r_n(\cos \theta_n+{ I}\sin \theta_n) \, (r_n<1,\, n\in \N^*)$, then ${\bf{\Phi}}_n, (n\in \N^*) $ is a slice regular  orthonormal system in  $H^2(\mathbb{B})$.

{\bf Proof.}
Recall the definition of the  inner product on the space $H^2(\mathbb{B})$: if
 $f(q) =\sum_{n\ge 0} q^na_n$ and  $g(q) =\sum_{n\ge 0} q^nb_n$, then their inner product is
$$\langle f, g\rangle=\sum_{n \ge 0}\overline{b_n}a_n=\lim_{r\to 1}\frac{1}{2\pi}\int_0^{2\pi}\overline{g(re^{I\theta})}f(re^{I\theta})d\theta $$ for any $I \in \mathbb{S}$.
Let  $I$ be the direction fixed by $a_n=r_ne^{\theta_n {\bf I}}=r_n(\cos \theta_n+{\bf I}\sin \theta_n)$. On $\mathbb{B}_I=\mathbb{B}\cap L_I$ the regular Blaschke product is  slice preserving i.e. $$\mathcal{B}_{a}(\overline{\mathbb{B}_I})\subset \overline{\mathbb{B}_I}.$$ 
Moreover, an easy computation shows that for $q\in \overline{\mathbb{B}_I}$: $$(1-e^{I\theta}\overline{a})^{-\ast}=(1-e^{I\theta}\overline{a})^{-1},\,\,\mathcal{B}_{a}(q)={B}_{a}(q).$$ From the  slice preserving property and the Splitting Lemma of slice regular functions follows that for  every $z=x+Iy\in \overline{\mathbb{B}_I}$ we have $\mathcal{B}_{a}(z)=F(z)$, where $F(z)$ is holomorphic in $\overline{\mathbb{B}_I}$. This implies that the slice regular $\ast$-product on this slice is equal to  the point-wise product of the factors, moreover on the  slice  $ \overline{\mathbb{B}_I}$ we have:
$$
{\bf{\Phi}}_n(z)={\sqrt{1-|a_n|^2}}\left (\ast\prod _{k=1}^{n-1}{\mathcal{B}}_{a_k}(z)\right )\ast (1-z\overline{a_n})^{-\ast}={\sqrt{1-|a_n|^2}}\left (\prod _{k=1}^{n-1}B_{a_k}(z)\right ) (1-z\overline{a_n})^{-1},
$$
and the order of the factors, can be  interchanged, because the pointwise product is commutative on the slice $\overline{\mathbb{B}_I}$.

Slice regular Blaschke functions and  classical Blaschke functions maps the unit ball into the unit ball, and the boundary into itself, consequently $|{B}_{a}(e^{I\theta})\overline{{B}_{a}(e^{I\theta})}|=1$. Taking into account these properties, the commutativity of the product on the slice $\overline{\mathbb{B}_I}$, and the Cauchy formula,  in the proof we can  follow the same line as in the complex case:
$$\langle {\bf{\Phi}}_n, {\bf{\Phi}}_n \rangle=\lim_{r\to 1_-}\frac{1}{2\pi}\int_0^{2\pi}\overline{{\bf{\Phi}}_n(re^{I\theta})}{\bf{\Phi}}_n(re^{I\theta})d\theta =$$ $$
{(1-|a_n|^2)} \frac{1}{2\pi}\int_0^{2\pi}\overline{(1-e^{I\theta}\overline{a_n})^{-1}}(1-e^{I\theta}\overline{a_n})^{-1}d\theta=1.$$
For $m>n$ we have:
$$\langle {\bf{\Phi}}_n, {\bf{\Phi}}_m \rangle=\lim_{r\to 1_-}\frac{1}{2\pi}\int_0^{2\pi}\overline{{\bf{\Phi}}_n(re^{I\theta})}{\bf{{\Phi}}}_m(re^{I\theta})d\theta =$$ 
$$
=\frac{1}{2\pi}\int_0^{2\pi}\overline{{\bf{\Phi}}_n(e^{I\theta})}{\bf{\Phi}}_m(e^{I\theta})d\theta=
$$
$$
=\frac{1}{2\pi}\int_0^{2\pi}\overline{{\sqrt{1-|a_n|^2}}\left (\prod _{k=1}^{n-1}B_{a_k }(e^{I\theta})\right ) (1-e^{I\theta}\overline{a_n})^{-1}}{\sqrt{1-|a_m|^2}}\left (\prod _{k=1}^{m-1}B_{a_k}(e^{I\theta})\right ) (1-e^{I\theta}\overline{a_m})^{-1}d\theta
$$
$$
={\sqrt{1-|a_n|^2}}{\sqrt{1-|a_m|^2}}\prod _{k=n}^{m-1}B_{a_k}(a_n)(1-a_n\overline{a_m})^{-1}=0.
$$

{\bf Theorem 4.2.}
If $\sum_{n\ge 0}(1 - |a_n|) = +\infty$, then the system  $\bf{\Phi}_n, (n\in \N^*) $ is complete in $H^2(\mathbb{B})$.

{\bf Proof.}
 To prove that the system is complete in $H^2(\mathbb{B})$ we need to prove the following implication: if for an  $f\in H^2(\mathbb{B})$ we have that $\langle f, {\bf{\Phi}}_n\rangle=0, \, n\in \mathbb{N}^*$, then $f\equiv 0$. According to the Splitting Lemma there exist two holomorphic functions $F,G : 
\mathbb{B}_I \to L_I$, such that for every
$z = x + yI \in \mathbb{B}_I$ , we have
$f_I (z) = F(z) + G(z)J$, where $J$ is orthogonal to $I$. Moreover $F, G \in H^2(\mathbb{B}_I)$. Then from $\langle f, \bf{\Phi}_n\rangle=0, \, n\in \mathbb{N}$ we get that $\langle F, {\bf{\Phi}}_n\rangle=0$ and $\langle G, {\bf{\Phi}}_n\rangle=0, \, n\in \mathbb{N}$. Because  on the slice  $\mathbb{B}_I$ the functions $F, G, \bf{\Phi}_n$ are slice preserving holomorphic functions,  analogue as in the case of the complex M-T, which is complete under the assumption of the theorem, we get that $F(z)=G(z)=0, \, z\in \mathbb{B}_I$. Consequently $f_I(z)=0, \, z\in \mathbb{B}_I$.

According to the Extension Lemma holomorphic functions defined on a domain 
$\Omega_I$, symmetric with respect to the real axis in the complex plane
$L_I$, extend uniquely to the smallest symmetric slice domain of $\mathbb{H}$ containing 
$\Omega_I$. We apply this to the unit ball $\mathbb{B}$ and his slice $\mathbb{B}_I$.
 Then for $f\in H^2(\mathbb{B})$ we have
$$
f(x + yJ) =\frac{1}{2}[f_I (x + yI) + f_I (x - yI)] + \frac{JI}{2}[f_I(x - yI) - f_I (x + yI)]$$
extends $f_I$  uniquely to a regular function $f : \mathbb{B} \to \mathbb{H}$. Taking into consideration that  $f_I(z)=0, \, z\in \mathbb{B}_I$, we have  $f(z)=0, \, z\in \mathbb{B}$.

\section{{ The properties of the projection operator  }}

\vskip5mm  

For $f\in H^2(\mathbb{B})$   according to the Splitting Lemma,  there exist two holomorphic functions $F,G : 
\mathbb{B}_I \to L_I$ , such that for every
$z = x + yI \in \mathbb{B}_I$ , we have
$f_I (z) = F(z) + G(z)J$, where $J$ is orthogonal to $I$. Moreover $F, G \in H^2(\mathbb{B}_I)$.
Let us consider the
boundary limit of functions $f$ in $H^2(\mathbb{B})$.    Similarly to the complex case,   for almost every $\theta \in  [0,2\pi)$, the limit
$$\lim_{r\to 1_-}f(re^{I\theta}) =  {f}_I (e^{I\theta})={F} (e^{I\theta})+ {G} (e^{I\theta})J$$
exists for all $I \in \mathbb{S}$ and in this case $ f_I(e^{I\theta}), {F} (e^{I\theta}), {G} (e^{I\theta})$ belong to $L^2(\partial \mathbb{B}_I )$.
 
 Let us consider the
orthogonal projection operator of an arbitrary function $f\in
H^2( \mathbb{B})$ on the subspace $V_n$ spanned by   the functions $\{{\bf \Phi}_k,\, k=1, \cdots, n\}$
\begin{equation}
\label{Pnf}
P_nf(z)=\sum_{k=1}^n {\bf \Phi}_{k}(z) \langle f, {\bf\Phi}_{k}\rangle , 
\end{equation}
where the value of the scalar product $\langle f, {\bf\Phi}_{k}\rangle$ is 
$$
\langle f, {\bf\Phi}_{k}\rangle=\lim_{r\to 1_-}\frac{1}{2\pi}\int_0^{2\pi}\overline{{\bf{\Phi}}_k(re^{I\theta})}f(re^{I\theta})d\theta=$$
$$\lim_{r\to 1_-}\frac{1}{2\pi}\int_0^{2\pi}\overline{{\bf{\Phi}}_k(re^{I\theta})}F(re^{I\theta})d\theta+ \lim_{r\to 1_-}\frac{1}{2\pi}\int_0^{2\pi}\overline{{\bf{\Phi}}_k(re^{I\theta})}G(re^{I\theta})d\theta J.
$$ 
On the  slice  $ \overline{\mathbb{B}_I}$ we have:
$$
{\bf{\Phi}}_n(z)=\frac1{\sqrt{1-|a_n|^2}}\left (\ast\prod _{k=1}^{n-1}{\mathcal{B}}_{a_k}(z)\right )\ast (1-z\overline{a_n})^{-\ast}=$$ $$\frac1{\sqrt{1-|a_n|^2}}\left (\prod _{k=1}^{n-1}B_{a_k}(z)\right ) (1-z\overline{a_n})^{-1}={{\Phi}}_n(z),
$$
consequently the coefficients of the projection operator can be expressed by
$$\langle f, {\bf\Phi}_{k}\rangle=\frac1{2\pi}\int_0^{2\pi}\overline{{{\Phi}}_k(e^{I\theta})}F(e^{I\theta})d\theta+ \frac{1}{2\pi}\int_0^{2\pi}\overline{{{\Phi}}_k(e^{I\theta})}G(e^{I\theta})d\theta J.
$$ 

If $\sum_{n\ge 0}(1 - |a_n|) = +\infty$, then the system  ${\bf{\Phi}}_n, (n\in \N^*) $ is complete in $H^2(\mathbb{B})$, this   implies that  for every $f\in H^2(\mathbb{B})$ the projection of $f$ on $V_n$ converges in norm to $f$, i.e. we have  
$$
\| f-P_nf \|\to 0, \quad n\to \infty. 
$$
 Since convergence in
norm implies uniform convergence inside the unit ball $\mathbb{B}$ on every compact subset, we
conclude that $ P_nf(z)\to  f(z)$ uniformly on every compact subset
of the unit ball.

{\bf Theorem 5.1.} { If the parameters of the slice regular Malmquist -Takenaka system are on the same slice, i.e., there exists $I\in \mathbb{S}$ such that $a_n=r_ne^{\theta_n { I}}\, (r_n<1,\, n\in \N^*)$, then for all $f\in H^2(\mathbb{B})$ the restriction of the projection operator $P_nf$  to  the slice $\mathbb{B}_I$ of the unit
ball} is an
interpolation operator  in the points $a_{\ell}=r_{\ell}e^{\theta_{\ell} { I}}\, ( {\ell}\in \{1, \cdots, n\})$.

{\bf Proof.} 
The restriction of the projection $P_nf$ to the slice $ \overline{\mathbb{B}_I}$ can be written in closed form as follows:
$$(P_n)_If(z)=\sum_{k=0}^n ({\bf \Phi}_{k})_I(z) \langle f, {\bf\Phi}_{k}\rangle =$$ $$\frac{1}{2\pi}\int_0^{2\pi} \sum_{k=0}^n ({ \Phi}_{k})_I(z)\overline{{{\Phi}}_k(e^{I\theta})}F(e^{I\theta})d\theta + \frac{1}{2\pi}\int_0^{2\pi} \sum_{k=0}^n ({ \Phi}_{k})_I(z)\overline{{{\Phi}}_k(e^{I\theta})}G(e^{I\theta})d\theta J. $$

For the Dirichlet kernel   of the classical extension over the set of the quaternions of the M-T system   we have an analogue of Darboux-Christoffel formula (see \cite{pa04}): 
$$
D^*_n(z,w):=\sum_{\ell=1}^n \Phi_{\ell}(z)(1-z\overline{w})\overline{\Phi_{\ell}(w)}=1-\prod_{\ell=1}^nB_{a_{\ell}}(z)\prod_{\ell=1}^N\overline{B_{a_{n-\ell+1}}(w)}, 
$$
for all $z,w\in {\overline{\mathbb B}},\,\, z\ne w$. 
Taking the restriction of the Dirichlet kernel to the slice $ {\overline{\mathbb {B_I}}}$, where the product is commutative, and using the slice preserving property of $\Phi_{\ell}(z)$ on ${\overline{\mathbb {B_I}}}$ 
we get that the restriction of the projection operator on the slice ${\overline{\mathbb {B_I}}}$ can be expressed very similar to the complex case:
$$(P_n)_If(z)= $$
$$
=\frac{1}{2\pi}\int_0^{2\pi}(1-ze^{-I\theta})^{-1}\left(1-\prod_{\ell=1}^nB_{a_{\ell}}(z)\prod_{\ell=1}^n\overline{B_{a_{n-\ell+1}}(e^{I\theta})}\right) F(e^{I\theta})d\theta +$$ $$+ \frac{1}{2\pi}\int_0^{2\pi} (1-ze^{-I\theta})^{-1}\left(1-\prod_{\ell=1}^nB_{a_{\ell}}(z)\prod_{\ell=1}^n\overline{B_{a_{n-\ell+1}}(e^{I\theta})}\right) G(e^{I\theta})d\theta J.
$$
From here using the Cauchy formula we get that the restriction of the projection $P_nf$ to the slice ${\overline{\mathbb {B_I}}}$ is an interpolation operator for the points $z=a_{\ell}, \, \ell=1, \cdots, n$. Indeed we have:
$$(P_n)_If(a_{\ell})=\frac{1}{2\pi}\int_0^{2\pi}(1-a_{\ell}e^{-I\theta})^{-1}F(e^{I\theta})d\theta + \frac{1}{2\pi}\int_0^{2\pi} (1-a_{\ell}e^{-I\theta})^{-1}G(e^{I\theta})d\theta J=$$
$$
F(a_{\ell})+G(a_{\ell}) J=f_I(a_{\ell}).
$$
{\bf Summary.} We have introduced   the slice regular analogue of the Malmquist-Takenaka system, and we proved that, similar to the complex case, this is a complete orthonormal system in the slice regular   Hardy space of the unit ball. 
  We proved that  the associated  projection  operator $(P_nf, \, n\in \Bbb N)$  
  is convergent in $H^2(\mathbb{B})$ norm  to $f$, and $ P_nf(z)\to  f(z)$
  uniformly on every compact subset of the unit ball.
   In the same time the restriction of $ P_nf(z)$ to a slice $\mathbb{B}_I$ is an  interpolation operator  of $f$.

\end{document}